\documentclass{amsart}
\usepackage[latin1]{inputenc}
\usepackage{amsmath}
\usepackage{amsfonts}
\usepackage{amssymb}
\usepackage{amsthm}
\usepackage{wrapfig}
\usepackage{url}

\title{Permutation Designs and Sequencing Highly Transitive Group Actions}
\author{Tad White}
\email{tad [preposition] super [delimiter] org}
\address{IDA Center for Computing Sciences, Bowie, MD 20715}

\newcommand{\Z}{\mathbf{Z}}

\newcommand{\Sym}[1]{\mathrm{Sym}(#1)}
\newcommand{\Alt}[1]{\mathrm{Alt}(#1)}
\newcommand{\Aff}[1]{\mathrm{Aff}(#1)}
\newcommand{\pgl}{\mathrm{PGL}_{2}}
\newcommand{\bx}{\mathbf{x}}
\newcommand{\br}{\mathbf{r}}

\newcommand{\carth}{Carthaginian \relax}

\newcommand{\fp}[2]{\relax{#1}^{\underline{#2}}}

\newcommand{\gf}[1]{\mathbb{F}_{#1}}
\newcommand{\PP}{\mathbb{P}}

\theoremstyle{plain}
\newtheorem{theorem}{Theorem}
\newtheorem{lemma}[theorem]{Lemma}

\newtheorem{conj}[theorem]{Conjecture}

\theoremstyle{definition}
\newtheorem{defn}[theorem]{Definition}
\let\df\textit

\long\def\comment#1{\relax}

\begin{document}


\begin{abstract}
We consider an experimental design problem for permutations: given a fixed set $X$, and an integer $t$,
construct a list $L$ of permutations of $X$ such that every ordered $t$-tuple of distinct elements of $X$
occurs as a consecutive subsequence of exactly one permutation in $L$. In this paper we focus on
solutions based on sharply transitive group actions, in effect generalizing Gordon's notion of group sequencing.
We give an explicit construction when $|X|$ is prime for the case $t=3$, and analyze
a branching algorithm for the general case which produces, for example, a rare design with $t=6$ based
on the Mathieu group $M_{12}$, and suggests that every sharply transitive group action leads to a solution,
apart from an explicit list of counterexamples. We state a number of conjectures and indicate
directions for future work.
\end{abstract}

\maketitle

\section{Introduction}
We consider here a problem of testing a sequential process in which a known set of $n$ events (or inputs) occur in some order, and we suspect that certain contiguous sequences of inputs may yield interesting behavior. Had we but world enough, and time, we would test all $n!$ permutations of the inputs. However, that requires time superexponential in $n$, so we settle for testing $t$-th order effects. Specifically, we can choose a small integer $t\ll n$, and try to construct a list of permutations which together contain each ordered $t$-tuple of distinct inputs exactly once as a contiguous subsequence. Numerology gives us hope for success: there are exactly
\[ \fp{n}{t} = n(n-1)\cdots(n-t+1)\]
possible $t$-tuples, and each $n$-permutation contains exactly $(n-t+1)$ of them.
Since $\fp{n}{t}/(n-t+1)=\fp{n}{t-1}$ is an integer, it is \textit{a priori} conceivable that this problem may
have a solution for any $n$ and $t$. This motivates the following definition:

\begin{defn} An \df{$(n,t)$-permutation design}, or \df{$(n,t)$-pd}, is a set of $\fp{n}{t-1}$ permutations of
an $n$-set, such that each ordered $t$-tuple occurs exactly once as a contiguous subsequence of one of the permutations.
\end{defn}
It is convenient to specify such a design as a $\fp{n}{t-1}\times n$ matrix whose rows are the permutations,
with the understanding that the order of the rows doesn't matter. Table \ref{table:aff5} displays a $(5,3)$-pd;
it consists of $20$ permutations of $\{1,\ldots,5\}$ such that each ordered triple occurs consecutively in exactly
one permutation. It was constructed from the affine permutations of the finite field $\gf5$; we will say much more
about the construction later. This design can be partitioned into a complete set of four mutually orthogonal Latin
squares, as shown. The dashed lines will be explained in section \ref{sec:nonsharp}.)

Observe that an $(n,1)$-pd is just a single permutation, and an $(n,n)$-pd is just an enumeration of the
symmetric group $\Sym n$. More interestingly, an $(n,2)$-pd is a set of $n$ permutations of an
$n$-set having the property that each ordered pair of distinct symbols occurs exactly once as a contiguous pair
somewhere in the set. Such objects are sometimes referred to as
\df{Tuscan squares} (\cite[IV.48]{colbourn2010crc}, \cite[p.~ 80]{denes91}.) For example, the set of rows of
the following matrix forms a $(7,2)$-pd (that is, a $7\times7$ Tuscan square) which is not a Latin square.
\begin{equation}
\begin{matrix}
1&2&3&4&5&6&7\\
2&5&7&4&6&1&3\\
3&7&2&1&6&5&4\\
4&1&7&6&3&5&2\\
5&3&6&2&4&7&1\\
6&4&2&7&3&1&5\\
7&5&1&4&3&2&6\\
\end{matrix}
\label{eq:tuscanexample}
\end{equation}
There is an extensive literature on Tuscan squares, so aside from discussing how these objects fit into
the present framework, this paper will focus on the cases $2<t<n$.

\newcommand{\foo}{\makebox[1ex]{\rule[-2pt]{0.1pt}{5pt}}\hspace{-1ex}\makebox[1ex]{\rule[5pt]{0.1pt}{4pt}}}
\begin{wraptable}[23]{r}{13ex}
\begin{center}
\begin{tabular}{|c|}
\hline
0 1 4\foo 2 3 \\ 1 2 0\foo 3 4 \\ 2 3 1\foo 4 0 \\ 3 4 2\foo 0 1\\ 4 0 3\foo 1 2 \\ \hline
0 2 3 4 1  \\ 1 3 4 0 2 \\ 2 4 0 1 3\\ 3 0 1 2 4\\ 4 1 2 3 0\\ \hline
0 3 2 1 4  \\ 1 4 3 2 0\\ 2 0 4 3 1\\ 3 1 0 4 2\\ 4 2 1 0 3 \\ \hline
0 4 1\foo3 2 \\ 1 0 2\foo 4 3\\ 2 1 3\foo 0 4\\ 3 2 4\foo 1 0 \\ 4 3 0\foo 2 1 \\ \hline
\end{tabular}
\parbox{1in}{\textsc{\ \strut\ Table \ref{table:aff5}.}}
\parbox{1in}{\ A $(5,3)$-pd}
\label{table:aff5}
\end{center}
\end{wraptable}
\stepcounter{table} 


In sections \ref{sec:groupsequencing} and \ref{sec:actionsequencing} we review group sequencings,
introduce sequencings of sharply transitive group actions, and use these to construct permutation
designs. In section \ref{sec:specificsequencings} we 
consider some specific families of actions; the subsequent sections discuss brute-force search, a duality
analogous to that between sequencings and terraces, and non-sharply transitive group actions.

The present work was inspired by discussions with David Cohen on his work (e.g.~\cite{cohen1998techniques})
on software testing designs. The most closely related work of which the author is aware is a paper of Chee et al \cite{chee2013sequence}, which considers a related
situation in which the $t$-tuples consist of distinct elements but need not occur contiguously.
The author is indebted to Katie Ahrens and Xander Faber for helpful discussions, insights and patience.

\section{Group sequencings}
\label{sec:groupsequencing}
We begin by reviewing the notion of a group sequencing, introduced by Basil Gordon in \cite{gordon61}.
Standard references include \cite[Ch.~3]{denes91}, \cite[IV.5.3]{colbourn2010crc},
\cite[Ch.~3]{chung201950} and the dynamic survey \cite{ollis13}.
Recall that a \df{Latin square} is an array in which each row \textbf{and} each column is a permutation
of a fixed set. A square which is simultaneously Latin and Tuscan is also called a \df{row-complete} Latin
square. Gordon introduced the notion of a group sequencing in order to construct such squares.
We review Gordon's fundamental construction using the now-standard terminology.

A \df{sequencing} of a group $G$ is an enumeration $(b_1,\ldots,b_n)$ of the elements of
$G$ such that the partial products $a_i$, with $a_i=b_1\cdots b_i$, are distinct. Note that this implies
$b_1=e$, the identity element of $G$. The sequence $(a_1,\ldots,a_n)$ is called a
\df{basic directed terrace} for $G$. The group $G$ is called \df{sequenceable} if there exists a
sequencing of $G$.

\begin{theorem}(Gordon, 1961) \cite{gordon61} If $(b_1,\ldots,b_n)$ is a sequencing of a group $G$ with associated basic directed terrace $\{a_i\}$, then the square matrix whose $(i,j)$ entry is $a_i^{-1}a_j$ is a row-complete (in fact, also column-complete) Latin square.
\end{theorem}

Observe that the entries of the Latin square are most simply described using the $a_i$'s rather than the $b_i$'s, and the $a_i$'s also contain each element of $G$ exactly once, so one might imagine a parallel universe in which the sequence $\{a_i\}$ receives first billing.
We will come back to this point in section \ref{sec:rhospace}.

\section{Sequencing sharply transitive group actions}
\label{sec:actionsequencing}
The basic observation is that \textit{if} we have a sharply highly transitive action of a group on a finite set,
and \textit{if} we can find a sequence enumerating the set appropriately, then the orbit of that sequence
under the group action yields a permutation design. We now make this observation precise.

Let $X$ be a set with an action of a group $G$ (i.e., a $G$-set). $G$ is said to act \df{transitively} on
$X$ if, given any $x, y\in X$, there exists $g\in G$ such that $gx=y$. $G$ acts coordinatewise on the
$k$-th Cartesian power $X^k$ of $X$; this action restricts to an action on the set $X^{(k)}\subset X^k$
 consisting of ordered $k$-tuples of \textit{distinct} elements of $X$. $G$ is said to act \df{$k$-transitively}
 on $X$
 if the action on $X^{(k)}$ is transitive; that is, given any sequences $\mathbf{x}=(x_1,\ldots,x_k)$ and
 $\mathbf{y}=(y_1,\ldots,y_k)$ of distinct elements of $X$, there exists some $g\in G$ such that
 $g x_i = y_i$ for all $i=1,\ldots,k$. The action of $G$ on $X$ is said to be \df{sharply $k$-transitive}
 \cite[V.1]{beth1999design} if this
 $g$ is unique; equivalently, if the stabilizer of $\mathbf{x}\in X^{(k)}$ consists solely of the identity element
 of $G$. We will refer to $X$ as a \df{(sharply) $k$-transitive $G$-set} if the action of $G$ is (sharply) $k$-transitive.

\subsubsection*{Sequencing of group actions.}
Let $X$ be a sharply $k$-transitive $G$-set with $n$ elements. A \df{sequencing} of $X$ is an enumeration $(x_1,\ldots,x_n)$ of the elements of $X$ such that each of the $n-k$ sequences
\begin{equation}
(x_1,\ldots,x_{k+1}),\ (x_2,\ldots,x_{k+2}),\ \ldots,\ (x_{n-k},\ldots,x_n) \tag{*}
\end{equation}
lie in distinct orbits of the action of $G$ on $X^{(k+1)}$.

Observe that if $G$ acts sharply $k$-transitively on $X$, then since the stabilizer of any $k$-tuple of distinct elements is trivial, $G$ has exactly $n-k$ orbits on $X^{(k+1)}$. So, since the $n-k$ $(k+1)$-tuples in (*) lie in distinct orbits, they in fact form a complete set of orbit representatives for the action of $G$ on $X^{(k+1)}$.

The following observation is now straightforward:
\begin{theorem}
\label{thm:gsetsequencing}
Let $\bx=(x_1,\ldots,x_n)$ be a sequencing of a sharply $k$-transitive $G$-set $X$. Then the orbit of $\bx$ under $G$ forms an $(n,k+1)$-permutation design.
\end{theorem}

\proof Let $T$ be a $(k+1)$-tuple in $X^{(k+1)}$. By assumption, the $n-k$ 
$(k+1)$-tuples in $(*)$ form a complete set of orbit representatives on $X^{(k+1)}$, so there
is a unique $(k+1)$-tuple in (*) such that some $g\in G$ sends it to $T$. By sharpness, that $g$
is unique, so $T$ occurs precisely once as a subsequence of an image of $\bx$ under $G$.\qed\\

One unavoidable downside of this construction is the dearth of sharply transitive actions, as quantified by the following classical result of Jordan:
\begin{theorem}[{\cite[Theorem~7.6A]{dixon2012permutation},\cite[Theorem~V.1.9]{beth1999design}}]
\label{thm:jordan}
The only sharply $k$-transitive permutation groups with $k\ge 4$ are the symmetric groups $\Sym{n}$ with $n\ge4$;
the alternating groups $\Alt n$ for $n\ge6$; and the Mathieu groups $M_{11}$ and $M_{12}$, which are $4$- and $5$- transitive, respectively.
\end{theorem}

We will describe completely the sequencing of the symmetric and alternating group actions, though these actions aren't particularly useful for constructing compact permutation designs. We will generally assume $k\ge1$, since taking $k=1$ yields precisely the well-studied theory of group sequencing. Thus, the only possible infinite families of interest must have $k=2$ or $k=3$. The simplest families of these types are the affine groups and the fractional linear groups respectively; we will discuss these particular groups in some detail.

\subsubsection*{\carth rectangles.}
Just as the Tuscan squares produced by Gordon's
construction have additional structure (in particular, they are always Latin squares), the $(n,k+1)$-permutation designs constructed via Theorem \ref{thm:gsetsequencing} have the following structure, which reduces to the Latin square
property when $k=1$.

We call a $\fp{n}{k}\times n$ matrix a $t$-\df{\carth rectangle} if its rows are permutations of a fixed set $X$ of $n$ symbols, and each $t$-tuple of columns together contains each of the $\fp{n}{k}$ elements of $X^{(k)}$
exactly once.%
\footnote{These designs are like Latin squares but are much ``longer.' Our terminology is inspired by the
longest surviving ancient Latin poem: Silius Italicus' \textit{Punica} \cite{enwiki:1026824677}, which recounts
the Second Punic War between Rome and Carthage.}
(It follows that the rows must be distinct.) It is immediate that if a group $G$ acts
sharply $k$-transitively on a set $X$, then the orbit of $X$ under $G$ (viewed as the $|G|\times|X|$ matrix whose $(g,x)$ entry is $gx$) forms a $k$-\carth rectangle (regardless of whether the action can be sequenced.)

There are standard criteria, such as the quadrangle condition \cite[\S2.1]{evans2018orthogonal}, for determining whether
a Latin square is constructed from a group in this way. The following condition reduces to the quadrangle condition
when $k=1$: we say that a matrix $R$ \df{satisfies the $k$-rectangle condition} provided, whenever two
$2\times (k+1)$ submatrices $A$ and $B$ of $R$ agree in $2k+1$ entries, the remaining entry agrees as well.
We defer the proof of the following result to appendix \ref{sec:carthageproof}.

\begin{theorem}
\label{thm:carthage}
A $k$-\carth rectangle $R$ arises from a sharply $k$-transitive group action if and only if $R$ satisfies the
$k$-rectangle condition. 
\end{theorem}

\subsubsection*{Equivalence and the reference sequence.}
We say that two sequencings $(x_i)$ and $(x'_i)$ of a $G$-set $X$ are \df{equivalent} if they are in the same
orbit under the diagonal action of $G$. We can identify a ``canonical'' element of the equivalence class
as follows.
Fix a \df{reference sequence} $\br=(r_1,\ldots,r_k)$ of distinct elements of $X$, and let $R=\{r_i\}\subset X$
be the set of reference elements. Sharp $k$-transitivity guarantees that any sequencing of $x$ is uniquely
equivalent to a unique sequencing that starts with the reference sequence. We say such a sequencing is in
\df{standard form}; we will generally assume all sequencings are in standard form. Thus, if we are looking
for a sequencing $x_1,\ldots,x_n$, we may as well think of $x_1,\ldots,x_k$ as being determined, and the
problem of finding a sequencing is really to find a suitable enumeration $x_{k+1},\ldots,x_n$ of $X\setminus R$.

\subsection{Orbit invariants}
\label{sec:orbits}

By a \df{complete orbit invariant} of an action of $G$,
we mean a function $\rho$ which is invariant under the action, and which takes different values on different orbits.
(That is, any other orbit invariant $\sigma$ factors as $\sigma=\pi\circ\rho$, for some $\pi$.)
We will assume that $G$ acts sharply $k$-transitively on an $n$-set $X$, so $X^{(k+1)}$ splits up into exactly
$n-k$ orbits under $G$. We will make frequent reference to the following general construction of a
complete orbit invariant 
$\rho$ on $X^{(k+1)}$. (When $G$ is the affine group, $\rho$ will turn out to be a ratio; when $G=\pgl(\gf q)$,
$\rho$ will be a cross-ratio.)

\subsubsection*{The standard orbit invariant.}
Given any $(k+1)$-tuple $\bx=(x_1,\ldots,x_{k+1})\in X^{(k+1)}$, we can use sharp $k$-transitivity to
map the first $k$ values to the reference sequence, and then define the \df{standard orbit invariant} to be the place
where $x_{k+1}$ lands. That is, we define
\begin{equation}
\label{def:rho}
\rho(\bx) = g x_{k+1},
\end{equation}
where $g\in G$ is the unique element sending $(x_1,\ldots,x_k)$ to $(r_1,\ldots,r_k)$. Note that $\rho(\bx)$ belongs to $X\setminus R$; each of the $n-k$ orbits maps to a different element of $X\setminus R$, so $\rho$ is a complete orbit invariant.

Effective computation of $\rho$ depends on computing $g$. If $G$ is a linear group, for example, $g$ is
typically straightforward to compute. If $G$ is handed to us as a list of permutations, we can sort or index
$G$ on $g^{-1}(r_1,\ldots,r_k)$ so that the correct $g$ can be computed with a single table look-up
with $(x_1,\ldots,x_k)$ as the key.

A sequencing of the $G$-set $X$ in standard form thus consists of a sequence of $n-k$ choices of
non-reference elements of $X$, such that the $n-k$ orbit representatives
\[ \rho_{1}:=\rho(x_1,\ldots,x_{k+1}),\ldots, \rho_{n-k}:=\rho(x_{n-k},\ldots,x_n) \]
are all distinct (and thus necessarily range over $X\setminus R$.)

We remark that formula (\ref{def:rho}) makes sense even when $x_{k+1}\in\{x_1,\ldots,x_k\}$.
Thus $\rho$ extends naturally to $\rho:X^{(k)}\times X\to X$; the same is then true for \textit{any}
complete orbit invariant. The following fact is immediate, and will prove useful later:

\begin{lemma}
\label{thm:injective}
Let $\rho:X^{(k)}\times X\to X$ be a complete orbit invariant, extended as above.
Then, given $(x_1,\ldots,x_k)\in X^{(k)}$ and $y\in X$, we can always uniquely solve the equation
\[ \rho(x_1,\ldots,x_{k+1})=y\]
for $x_{k+1}$.
\end{lemma}

It is worth taking a moment to explicitly recover the notion of a group sequencing by considering the special
case $k=1$. A sharply $1$-transitive (left) action of $G$ on a set $X$ is equivalent, given a choice of basepoint,
to the regular (left) action of $G$ on itself by multiplication.
We can take the (one-long) reference sequence in $G$ to consist of the identity $e$
of $G$. The standard construction above then associates to a pair $(x_1,x_2)\in X^{(2)}$ the orbit invariant
$x_1^{-1}x_2$. Thus a sequencing of a $1$-transitive (i.e. regular) group action is an enumeration
$(x_1,\ldots,x_n)$ of $G$ such that the values $\rho_i=x_{i-1}^{-1}x_{i}$, for $i\ge 2$, are all distinct.
We can assume the sequencing is in standard form, i.e.~$x_1=e$.
(Note that none of the $\rho$'s can equal $e$, so we can augment the sequence of $\rho$'s if we like by defining $\rho_1=e$.) The augmented sequence $\{\rho_i\}$ of orbit invariants corresponds to what is called a group sequencing in the literature. Observe that $x_i= \rho_1\rho_2\cdots\rho_i$; these must all be distinct, so the $x_1^{-1}x_i$'s correspond to the standard notion of a basic terrace. If $x_1=e$ so that the sequencing is in standard form, then the
$x_i$'s form a basic terrace.

\section{Some specific group sequencings}
\label{sec:specificsequencings}
In this section we discuss sequencings of some infinite classes of sharply transitive group actions: the alternating groups, the affine groups, and the fractional linear groups.

\subsection{The alternating group}
\label{sec:alternating}
The group $G=\Alt n$, consisting of the even permutations, acts naturally on an $n$-set $X$. This action is sharply $(n-2)$-transitive, so we may ask whether a $(n,n-1)$-permutation design can be constructed by sequencing this action.

\begin{theorem}
\label{thm:alternating}
The natural action of $\Alt n$ on $X=\{1,\ldots,n\}$ can be sequenced precisely when $n$ is even. In this case, any enumeration of $X$ forms a sequencing.
\end{theorem}

\proof Any $(n-1)$-tuple $(x_1,\ldots,x_{n-1})$ has a unique extension to a permutation $(x_1,\ldots, x_n)$. The sign of this permutation serves as an orbit invariant for the action on $X^{(n-1)}$. Now suppose $(x_1,\ldots,x_n)$ is a sequencing of the action; this means $(x_1,\ldots,x_{n-1})$ and $(x_2,\ldots,x_n)$ lie in different orbits on $X^{(n-1)}$, which happens precisely when the two permutations $(x_1,\ldots, x_n)$ and $(x_2,\ldots,x_n,x_1)$ have different signs, which happens precisely when $n$ is even.

In particular, when $n$ is even, any permutation $\pi$ of $\{1,\ldots,n\}$ is a sequencing of the action of $\Alt n$. There are two equivalence classes of sequencings, determined by $\textrm{sign}(\pi)$; in particular, the set of even (or odd) permutations forms an $(n,n-1)$-permutation design when $n$ is even.\qed

When $n$ is odd, we cannot in general say anything about the existence of such a design; we can say only that this construction does not apply.

\subsection{The affine group}
\label{sec:affine}
We discuss the familiar case of the natural action of the affine group over a finite field $\gf q$ in some detail. We fix a prime power $q=p^r$, and take $G$ to be the group $\Aff{\gf q}$ of affine self-maps of $\gf q$; that is,
\[G=\{ x\mapsto ax+b\mid a\in\gf{q}^*, b\in\gf{q}\}.\]
$G$ acts sharply $2$-transitively on $X=\gf q$: given $x_1\ne x_2$ and $y_1\ne y_2$, there is a unique affine map
taking $(x_1,x_2)$ to $(y_1,y_2)$, namely
\[x\mapsto (y_2-y_1)/(x_2-x_1)\, (x-x_1) + y_1.\]
A sequencing of the $2$-transitive $G$-set $\gf q$ will yield a $(q,3)$-permutation design by Theorem \ref{thm:gsetsequencing}.

We start by selecting a convenient orbit invariant for the action of $G$ on $X^{(3)}$.
Referring to the general construction in section \ref{sec:orbits}, we take $(1,0)$ to be our standard reference points in $X=\gf q$. Then the unique map taking $(1,0)$ to $(x_1,x_2)$ is $g:z\mapsto (x_1-x_2)z+x_2$, so the recipe in section \ref{sec:orbits} gives
\[ \rho(x_1,x_2,x_3) = g^{-1} x_3 = (x_3-x_2) / (x_1 - x_2). \]
We will refer to this simply as ``the ratio'' determined by $(x_1,x_2,x_3)$. By the remarks in section \ref{sec:orbits}, if the $x_i$'s are distinct, this ratio lies in $\gf q\setminus\{1,0\}.$

The special case of a prime field is especially nice. 
\begin{theorem}
\label{thm:affineconstruction}
For any prime $p>2$, the sequence $(0,1,p-1,2,p-2,\ldots)$ of elements of $\gf p$ is a sequencing of the
$2$-transitive action of
$\Aff{\gf p}$ on $\gf p$. In particular, a $(p,3)$-permutation design exists for every odd prime $p$.
\end{theorem}

\proof
The successive differences $x_i-x_{i-1}$ are $(1,-2,3,-4,5,\ldots)$, so the ratios $\rho_i=(x_i-x_{i-1})/(x_{i-2}-x_{i-1})$ take the values $(2/1, 3/2, 4/3, 5/4, \ldots)$, with the division taking place in $\gf p$. But the map $x\mapsto x/(x-1)$ is injective for $x\ne 1$, so the $p-2$ ratios $x/(x-1)$ mod $p$, with $2\le x\le p-1$, are distinct.
\qed

Table \ref{table:aff5} illustrates this construction for $p=5$. The corresponding $(5,3)$-pd consists of $20$ permutations of $\{0,1,2,3,4\}$ given by the images of the sequence $(0,1,4,2,3)$ under the 20 affine functions $x\mapsto ax+b$, with $a\in\gf p^*$, $b\in\gf p$. If we group these maps by the multiplier $a$, what we see is a set of four mutually orthogonal Latin squares\footnote{A pair $(a_{ij})$ and $(b_{ij})$ of Latin squares with symbol set $X$ is said to be \df{orthogonal} iff the ordered pairs $(a_{ij},b_{ij})$ range over $X^2$; see, for example, \cite[Ch.~1]{denes91}.},
whose rows taken together contain each ordered $3$-tuple from $\gf p$ exactly once.

This design has the very useful property of being \df{indexable}: given a triple $(x_1,x_2,x_3)$ of
distinct elements of $\gf p$, it is easy to locate that triple in the table: compute the ratio (orbit invariant)
$\rho=(x_3-x_2)/(x_2-x_1)$ to determine whether the triple is on the left, in the middle, or on the right
of its row. The triple is an affine image of either $(0,1,4)$, $(1,4,2)$ or $(4,2,3)$ respectively; the
correct affine map can be found by solving two linear equations.

Unfortunately, the simple construction which works for prime fields doesn't seem to extend in an
obvious way to more general finite fields. We have constructed sequencings by direct search for
all fields of order at most $100$; those of order at most $50$ can be found in appendix \ref{sec:explicitsequencings}.

\subsection{The fractional linear group}
\label{sec:PGL2}
The natural action of $\pgl(\gf q)$ is a fruitful place to look for $(*,4)$-permutation designs.
Let $\gf q$ be a field of order $q=p^r$. $\mathrm{GL}_2(\gf q)$ denotes the group of $2\times 2$ invertible matrices
over $\gf q$; $\pgl(\gf q)$ is the quotient of $\mathrm{GL}_2(\gf q)$ by its center (i.e.~nonzero scalar multiples of the identity.) The order of $\pgl(\gf q)$ is $(q^2-1)(q^2-q)/(q-1) = \fp{(q+1)}3$
 (see, e.g., \cite[III.5.20.d]{beth1999design}).

The projective line $\PP^1(\gf q)$ is
$$\PP^1(\gf q)\ =\ (\gf{q}^2\setminus \{(0,0)\}) / (\bx \sim \lambda\bx) \ = \ \gf{q}\cup\{\infty\}$$
where we identify  $[(x,1)]$ with $x\in\gf q$ and $[(1,0)]$ with $\infty$.
$\pgl(\gf q)$ acts on $\PP^1(\gf q)$ via $[A] : [\bx] \mapsto [A\bx] $.

One typically writes $[(x,y)]$ as $x/y$, so that $A=\begin{pmatrix}a&b\\ c&d\end{pmatrix}$ sends
$ x/y $ to $(ax+by)/(cx+dy)$. That is, we can think of $\pgl(\gf q)$ as the group of fractional linear
transformations of $\PP^1(\gf q)$. In this view, the stabilizer of $\infty$ consists precisely of the affine
transformations of $F$; we already know that $\Aff{\gf q}$ acts sharply $2$-transitively on $\gf q$, so the
action of $\pgl(\gf q)$ on $\PP^1(\gf q)$ is sharply $3$-transitive. Accordingly, a sequencing of the action
will yield an $(n,4)$-permutation design, where $n=q+1$.

It is conventional to take $R=(0,1,\infty)$ as the reference sequence (as in section \ref{sec:orbits}.) The standard orbit invariant then takes the form of the following cross ratio
$$\rho(w,x,y,z) := \frac{(z-w)(y-x)}{(x-w)(y-z)}$$
since $\rho(0,1,\infty,z) = z$ for $z\notin R$. Thus a sequencing of this action consists of an enumeration $(x_1,\ldots,x_n)$ of $\PP^1(\gf q)$
for which the cross ratios $\rho(x_1,x_2,x_3,x_4),\ldots$, $\rho(x_{n-3},x_{n-2},x_{n-1},x_n)$ are distinct and range over $\PP^1(F)\setminus \{0,1,\infty\}$.

Unlike the case of the affine group actions, we do not know a general sequencing of the fractional linear group actions. Here is the state of affairs:
\begin{itemize}
\item $\pgl(\gf3)$ can be sequenced trivially since there is only one cross ratio to compute. That is, a sharply $3$-transitive action on a $4$-element set is in fact $4$-transitive, and any ordering of $\PP^1(\gf3)$ is a sequencing.
\item When $q=4, 5,$ or $7$, no sequencing exists. This can be seen by direct search or case-by-case analysis; however, we'll give an alternate argument in section \ref{sec:rhospace} for the $q=4$ and $q=5$ cases. (The case $q=4$ can also
be dispatched easily since it is just the natural action of $\Alt{5}$.)
\item We have found sequencings by direct search for all $\gf q$ with $8\le q<100$,  The discussion in section \ref{sec:search} suggests that sequencings should exist whenever $q\ge 8$. 
\end{itemize}

\section{Sequencing by depth-first search}
\label{sec:search}
If we want to sequence a group action for which we don't have a general construction, we can resort to direct search.

Fix a sharply $k$-transitive action of a group $G$ on an $n$-set $X$. We seek an enumeration $\mathbf{x}=(x_1,\ldots,x_n)$ of $X$ such that the orbit representatives $\rho_{1}=\rho(x_1,\ldots,x_{k+1})$, $\ldots$, $\rho_{n-k}=\rho(x_{n-k},\ldots,x_n)$ are distinct. As in section \ref{sec:orbits}, it is useful to fix a sequence of reference elements $\br=(r_1,\ldots,r_k)\in X^{(r)}$, and assume without loss that  $x_i=r_i$ for $i\le k$, so that the the standard orbit invariants relative to $\br$ of the $G$-action on $X^{(k+1)}$ correspond naturally to $X\setminus R$, where $R=\{r_i\}$.

Algorithmically we will use depth-first search on the search tree implied here. Initially, we have $m=n-k$ choices for $x_{k+1}$, each of which produce a valid $\rho_1$ with probability 1 since none of our $m$ $\rho$'s have been used yet. At the next step, we have $m-1$ choices for $x_{k+2}$; at this point one of the $m$ orbit representatives has been used, so each of the $m-1$ choices is valid with probability $(m-1)/m$. At the stage at which we have $m-i$ choices, each choice is valid with probability $(m-i)/m$, so the expected number of extensions at the $i$-th stage is $(m-i)^2/m$.

If all these extension counts were independent, we could multiply the expected values together to estimate the number of sequencings extending a given $k$-long prefix:
\begin{equation}
\label{eq:expseq}
 E(S)=\prod_{i=0}^{m-1}\frac{(m-i)^2}{m} = \frac{m!^2}{m^m}
\end{equation}
If we truncate the product at $j$, we get the expected number of nodes at depth $j$ in the search tree, so we can estimate the search tree size as
\begin{equation}
\label{eq:expnodes}
E(T) = \sum_{j=0}^m \frac{(\fp{m}{j})^2}{m^j}.
\end{equation}

Some skepticism is understandable, since the estimate not only doesn't depend on the group, it only depends on the number $m=n-k$ of free choices in the process. We remark that three of the natural actions we have described related to a finite field $\gf q$ -- the $1$-transitive action of $\gf{q}^*$ on itself, the $2$-transitive action of $\Aff{\gf q}$ on $F$, and the $3$-transitive action of $\pgl(\gf q)$ on $\PP^1(\gf q)$ -- have the same value of $m=q-2$, so based on this estimate, they should asymptotically exhibit similar behavior; this is something we can check empirically.

Table \ref{table:dfscounts} gives the results of experiments with small values of $m$, where we can run the DFS to completion; the actual counts are in the ballpark of the predicted values. One naturally conjectures that the natural action of $\Aff{\gf q}$ is always sequenceable, and the natural action of $\pgl(\gf q)$ is sequenceable for $q>7$. We have verified these conjectures by direct calculation for $q\le100$.%
\footnote{We also note that $\Alt n$ has either 0 or 2 sequencings; here $m=2$ so (\ref{eq:expseq}) predicts one sequencing. It's never correct, but only off by one, and exactly right on average!}

\begin{table}
$$\begin{tabular}{|c|c|c|c|c|c|c|c|c|}
\hline
$q$&\multicolumn{2}{|c|}{predicted}&\multicolumn{2}{|c|}{$\gf q^{*}$}&%
\multicolumn{2}{|c|}{$\Aff{\gf q}$}&\multicolumn{2}{|c|}{$\pgl(\gf q)$}\\
\hline
  &S&T&S&T&S&T&S&T\\
\hline
3& 1& 2& 1& 2& 1& 2& 1& 2\\
4& 1& 4& 0& 1& 2& 5& 0& 3\\
5& 1.33& 9.33& 2& 10& 1& 8& 0& 8\\
7& 4.6&  78.4& 4& 66& 12& 105& 0&72\\
8& 11.1& 276.4& 0& 271& 18& 301& 12& 271\\
9& 30.8& 1091& 24& 994& 52& 1094& 20& 928\\
11& 340& 22853& 288& 19250& 493& 23360& 318& 22164\\
13& 5585& 671052& 3856& 547746& 7374& 743906& 5660& 674890\\
\hline
\end{tabular}$$
\caption{Experimental results for computing sequencings. $S$ and $T$ are the number of sequencings in standard form and nodes in the search tree, respectively. The predicted values come from equations (\ref{eq:expseq}) and (\ref{eq:expnodes}); the remaining columns are the actual answers from a complete traversal of the DFS tree.
The zeros under the $\gf q^{*}$ action come from Gordon's theorem, since these are involution-free abelian groups.}
\label{table:dfscounts}
\end{table}

If the estimate (\ref{eq:expseq}) is even in the right ballpark, a sharply transitive action should be sequenceable provided $m=n-k$ is large. By Theorem \ref{thm:jordan}, apart from the symmetric and alternating groups, $k$ is universally
bounded, so we should find any non-sequenceable actions by looking at small values of $n$. Colbourn et al
\cite[VI.1.4]{colbourn2010crc} list all primitive group actions with $n\le 50$; we were able to sequence every
sharply $k$-transitive action in the table with $k\ge2$ except for the counterexamples already identified.
(We exhibit these sequencings in appendix \ref{sec:explicitsequencings}.) The only examples which are not sequenceable are the ones
we have already discussed. Accordingly, we make the following conjecture:

\begin{conj}
\label{conj:existence}
Every sharply $k$-transitive action of a group $G$, with $k\ge2$, can be sequenced, with the following exceptions:
\begin{itemize}
\item $G=\Alt{n}$, acting $(n-2)$-transitively on $\{1,\ldots,n\}$, for $n$ odd; and
\item $G=\pgl(\gf q)$, acting $3$-transitively on $\PP^1(\gf q)$, for $q\in\{4,5,7\}$.
\end{itemize}
\end{conj}

\section{Putting the orbit invariants first}
\label{sec:rhospace}

We have been viewing the problem of sequencing a sharply $k$-transitive $G$-set $X$ as picking an enumeration of $X$, and checking whether the orbit invariants $\{\rho_i\}$ on the $(k+1)$-long windows are distinct. We could also turn this around: given a sequence $\{\rho_i\}$ of distinct orbit invariants, can they arise as the values of $\rho$ on the windows of some sequence of distinct $x_i$'s. Of course distinctness is the only challenge here, as we've seen that $\rho$ is 1-1 when all arguments but the last are held fixed. The situation is particularly nice in the affine and fractional linear cases; as an application, we will give simple proofs that the natural actions of $\pgl(\gf4)$ and $\pgl(\gf5)$ cannot be sequenced.

\subsection{The affine case}
We will continue to use the ``standard'' orbit invariant $\rho(x_1,x_2,x_3)=(x_3-x_2)/(x_1-x_2)$ defined in section \ref{sec:affine}, though in principle any one will do. Given $\rho_1,\ldots,\rho_i$, we would like to know if there exist distinct $x_1,\ldots,x_{i+2}$ such that the $\rho_j$'s can be realized as the values of $\rho$ on the three-long windows of $x_1,\ldots,x_{i+2}$. As always, we may assume by $2$-transitivity that $x_1=1$ and $x_2=0$.

If $i=1$ this is easy: as long as $\rho_1\notin\{1,0\}$, we can (by Lemma \ref{thm:injective}) uniquely choose $x_3=\rho_1$ so that $\rho(1,0,x_3)=\rho_1$. So consider the case $i=2$. Given $\rho_1$ and $\rho_2$, and a putative sequence $x_1,\ldots, x_4$, we can ensure $x_1,x_2,x_3$ are distinct by checking $\rho_1$ as in the $i=1$ case, and that $x_2,x_3,x_4$ are distinct by checking $\rho_2$. So it only remains to determine whether $x_1=x_4$. If this does in fact hold, then
\begin{equation}
\label{eq:affinerhos}
\rho_1=(x_3-x_2)/(x_1-x_2)\textrm{\ and\ }\rho_2=(x_1-x_3)/(x_2-x_3).
\end{equation}
We can eliminate $x_1$, $x_2$, and $x_3$ from these two equations (using Gr\"obner bases, for example); we find that $1-\rho_1+\rho_1\rho_2=0$ (still assuming that $x_4=x_1$.)%
\footnote{%
It may seem a little mysterious that one can eliminate three variables from the two equations (\ref{eq:affinerhos}). This is just a restatement of the fact that $\rho(x_2,x_3,x_1)=1-1/\rho(x_1,x_2,x_3)$. The (probably better known) analogous fact in the fractional linear case is that the value of a cross-ratio, if one permutes the arguments, can be computed as a simple function of the original cross-ratio. We will say more about this mystery when discussing the case of general permutation groups.}
But it follows from Lemma \ref{thm:injective}, together with the assumption that $\rho_1\ne0$, that if $x_4\ne x_1$ then $1-\rho_1+\rho_1\rho_2\ne0$. So the vanishing of this polynomial precisely detects whether $x_1=x_4$.

\subsubsection*{Distinguishing polynomials}
We will refer to the polynomial $d_2(\rho_1,\rho_2) := 1 - \rho_1 + \rho_1\rho_2$ as a
\df{distinguishing polynomial}, since we can evaluate it on the orbit invariants of the two adjacent triples in a
$4$-long window
$(x_i,\ldots,x_{i+3})$ to determine whether the first and last values in the window ($x_i$ and $x_{i+3}$)
are distinct, assuming all other pairs of $x_j$'s in the window are distinct. Similarly, for any $s\ge1$, we observe that
\begin{equation}
\label{eq:affinemagic}
(x_1-x_{s+2})/(x_1-x_2) = 1 - \rho_1 + \rho_1\rho_2 -\cdots \pm\rho_1\rho_2\cdots\rho_{s} =: d_s(\rho_1,\ldots,\rho_{s}).
\end{equation}
The field characteristic or size does not enter into these computations, so the same polynomials $s$ apply to the affine action over any field. If we change our choice of orbit invariant (e.g. from $(x_3-x_2)/(x_1-x_2)$ to $(x_3-x_1)/(x_2-x_1)$, we may get a different sequence of polynomials $d_s$, but everything still works the same way.

We can use this, algorithmically, to search for a sequence of distinct values $\rho_1,\ldots,\rho_{n-2}$ of $\gf q\setminus\{1,0\}$ such that for each pair $i<j$,
\[ d_{j-i+1}(\rho_i,\rho_{i+1},\ldots,\rho_j)\ne 0.\]
This can be done reasonably efficiently, since the polynomials $d_s$ can be computed recursively:
\begin{equation}
\label{eq:recursivemagic}
d_s(x_1,\ldots,x_{s}) = 1-\rho_1 d_{s-1}(x_2,\ldots,x_{s}).
\end{equation}
Accordingly, each required value of $d_s$ can be computed in constant time given the values for smaller $s$, so the running time of the depth-first search is proportional to the number of nodes searched. Finally, once a suitable ordering $\{\rho_i\}$ has been found, Lemma \ref{thm:injective} guarantees that one can construct from it the corresponding sequence of $x_i$'s, which is unique up to equivalence.
 
\subsection{The fractional linear case}
This case works the same way, except that the orbit invariant $\rho$ is now a cross ratio, taking four arguments. The distinguishing polynomials $d_s$ can be computed in the same manner. For the specific cross ratio used in section
\ref{sec:PGL2}, we find the following:
$$
\begin{array}{ccccc}
s & d_s \textrm{(as a function of the $\rho_i$'s)} & d_s \textrm{(as a function of the $x_i$'s)}\\
\hline
 1 & \rho_1 & \frac{\left(x_3-x_2\right) \left(x_4-x_1\right)}{\left(x_2-x_1\right)
   \left(x_3-x_4\right)} \\
 2 & \rho_1 \rho_2-1 & \frac{\left(x_2-x_4\right) \left(x_1-x_5\right)}{\left(x_1-x_2\right)
   \left(x_4-x_5\right)} \\
 3 & \rho_1 \rho_2 \rho_3 -\rho_1-\rho_3+1 & \frac{\left(x_4-x_2\right)
   \left(x_3-x_5\right) \left(x_1-x_6\right)}{\left(x_1-x_2\right) \left(x_3-x_4\right)
   \left(x_5-x_6\right)} \\
 4 & \rho_1 \rho_2 \rho_3 \rho_4 - \rho_1\rho_2 -\rho_1\rho_4 -\rho_3 \rho_4 +
   \rho_1 + \rho_4 & \frac{\left(x_4-x_2\right) \left(x_3-x_5\right) \left(x_4-x_6\right)
   \left(x_1-x_7\right)}{\left(x_1-x_2\right) \left(x_4-x_3\right) \left(x_4-x_5\right)
   \left(x_6-x_7\right)} \\
\end{array}
$$
In general, it turns out that $d_s$ (when expanded) is a sum of monomials, each of which is obtained from $\rho_1\cdots\rho_{s}$ by striking out some intervals $I$ of consecutive variables, and whose coefficient is the product of $1-\left(\textrm{length}(I) \pmod3\right)$ over these intervals. We will omit the proof, since we won't be making use of this observation. We don't need to deal with the complexity of the expanded polynomials in practice, since
they can be computed recursively just as in the affine case.

To illustrate the utility of the distinguishing polynomials, we use them to give simple proofs that the natural
actions of $\pgl(\gf4)$ and $\pgl(\gf5)$ cannot be sequenced, without recourse to case-by-case analysis.

For $\pgl(\gf4)$, we seek two distinct cross ratios $\rho_1$ and $\rho_2$; we must have
\[ \{\rho_1,\rho_2\} = \PP(\gf4)\setminus\{0,1,\infty\} = \{\alpha,\alpha^2\}\]
where $\alpha\in\gf4$ is a root of $\alpha^3-1=0$. But this forces $\rho_1$ and $\rho_2$ to be reciprocals
of each other, contradicting the condition $d_2(\rho_1,\rho_2)=\rho_1\rho_2-1\ne0$.

Similarly, for $\pgl(\gf5)$, we seek a sequence $(\rho_1,\rho_2,\rho_3)$ of three distinct cross ratios in $\PP^1(\gf5)\setminus\{0,1,\infty\}$ on which the $d_s$'s don't vanish. The $\rho_i$'s must be $2$, $3$ and $4$ in some order. By nonvanishing of $d_2$, $2$ and $3$ cannot be consecutive (as they are reciprocals mod $5$), so $\{\rho_1,\rho_3\}=\{2,3\}$. But by Fermat's Little Theorem, $2\cdot 3\cdot 4\equiv -1\pmod5$, so $d_3(\rho_1,\rho_2,\rho_3)=\rho_1+\rho_3 = 2+3 = 0 \pmod5$. Thus $d_3$ is forced to vanish, and no sequencing can exist.

\subsection{General sharply transitive permutation groups}
One can ask whether, in the setting of a general sharply $k$-transitive $G$-set $X$, we can still determine whether a sequence of $\rho$'s can arise as the orbit invariants of $(k+1)$-long windows of a sequence of distinct elements of $X$. Of course in this setting we don't have polynomial rings floating around, but it turns out that for any given $i\ne j$, there is still an equation in the $\rho$'s that holds precisely when $x_i=x_j$.

The first key point is that if $\rho$ and $\rho'$ are both complete orbit invariants for the action on $X^{(k+1)}$
(with values in $X\setminus R$, without loss of generality, where $R$ is as in section \ref{sec:orbits}), then $\rho(x_1,\ldots,x_{k+1})$ can be computed only from the value $\rho'(x_1,\ldots,x_{k+1})$: you don't need to know anything else about the $x$'s since $\rho'$ completely determines the orbit. In particular, there exists some function $f:X\to X$ such that $\rho'=f(\rho)$.

We now observe that if $\rho:X^{(k+1)}\to X\setminus R$ is a complete orbit invariant, and $\sigma\in\mathrm{Sym}(k+1)$, then
$\rho(x_{\sigma(1)},\ldots,x_{\sigma(k+1)})$ is also a complete orbit invariant, for
\begin{align*}
&\quad\rho(x_{\sigma(1)},\ldots,x_{\sigma(k+1)})=\rho(y_{\sigma(1)},\ldots,y_{\sigma(k+1)})\\
\iff&\quad(x_{\sigma(1)},\ldots,x_{\sigma(k+1)}) = g(y_{\sigma(1)},\ldots,y_{\sigma(k+1)}) \textrm{\ for some $g\in G$}\\
\iff&\quad(x_1,\ldots,x_{k+1})= g(y_1,\ldots,y_{k+1}).
\end{align*}

Applying the preceding observation, there must be some (necessarily bijective) function $f:X\setminus R\to X\setminus R$ such that $(\sigma^*\rho)(\mathbf{x}) = \rho(\sigma\mathbf{x}) = f(\rho(\mathbf{x}))$. If we denote this $f$ by $f_\sigma$, we have
\[\sigma^*\rho = f_\sigma\circ\rho.\]
(This explains why, for example, permuting the arguments to a cross-ratio gives a value which can be computed as a function of the value of the original cross-ratio. If $\rho(w,x,y,z)=\frac{(z-w)(y-x)}{(x-w)(y-z)}$, and we take $\sigma$ to be a $4$-cycle permuting the variables, we have $\rho(x,y,z,w)=1/\rho(y,z,w,x)$, so that $f_\sigma$ is the function $x\mapsto 1/x$.)

Finally, we need the observation that a complete orbit invariant for the action of $G$ on $X^{(k+i)}$ is given by the sequence of windowed values $(\rho_1,\ldots,\rho_i)$.

We can now describe how to determine from the sequence $\rho_1,\ldots,\rho_i$ corresponding to the orbit invariants of $(k+1)$-long windows of a sequence $x_1,\ldots,x_{k+i}$, whether the $x_j$'s are distinct. Suppose we have already checked all pairs of elements for distinctness except for the pair $(x_1,x_{k+i})$. Following the argument in the affine case, suppose for the moment that $x_1=x_{k+i}$. Then the two $(k+i-1)$-long windows of $x$'s differ only by a cyclic shift $\sigma$. Accordingly, the sequences $(\rho_1,\ldots,\rho_{i-1})$ and $(\rho_2,\ldots, \rho_i)$ are both complete orbit invariants for the action of $G$ on $X^{(k+i-1)}$. That means the second must be a function $f_\sigma$ of the first; in particular $\rho_i$ can be expressed as $g(\rho_1,\ldots,\rho_{i-1})$ for some $g$. As before, $\rho_i$ is $1-1$ when its last argument varies, so the equation $\rho_i = g(\rho_1,\ldots,\rho_{i-1})$ holds precisely when $x_1=x_{k+i}$.
This equation generalizes, in principle, the vanishing of the distinguishing polynomials that arise in the affine and fractional linear cases.

\section{Non-sharply transitive actions}
\label{sec:nonsharp}
Sharply transitive actions are hard to come by.
If we are given an action of $G$ on $X$ which is transitive, but not sharply transitive, things become more complicated. We mention here a natural extension of the ideas which sometimes applies.

Recall that the action is sharply $k$-transitive if it is $k$-transitive and the pointwise stabilizer $G_K$ of any $k$-subset $K$ of $X$ is trivial. More generally, we define the action to be \df{flatly $k$-transitive} if it is $k$-transitive, and the orbits of $G_K$ on $X\setminus K$ all have the same size. (So if that size is $1$, the action is sharply $k$-transitive.) We will use the adjective phrase ``$s$-flatly $k$-transitive'' if we need to emphasize the common orbit size $s$.

Note that if the action is $s$-flatly transitive and $|X|=n$, then the order of $G$ is $s\fp{n}{k}$. There are (at least) two notions of group sequencing that make sense in this situation, both of which lead to analogues of a permutation design.

A \df{long sequencing} of a flatly $k$-transitive group action is an enumeration $(x_i)$ of $X$ such that the orbit invariants on the $(k+1)$-long windows $(x_i,\ldots,x_{k+i})$ take each value exactly $s$ times. It is easy to see that the orbit under $G$ of a long sequencing of $X$ forms a list of permutations of $X$ which together contain every $(k+1)$-long sequence of distinct elements exactly $s$ times; we'll call this a \df{wide} design.

A \df{short sequencing} of a flatly $k$-transitive group action is a sequence $(x_i)$ of $k+(n-k)/s$ distinct elements of $X$ such that the $(n-k)/s$ orbit invariants $\rho_i$ on the $(k+1)$-long windows are distinct (and thus each of the $(n-k)/s$ possible values of $\rho_i$ appears exactly once.) It is easy to see that the orbit under $G$ of a short sequencing of $X$ forms a list of partial permutations of $X$, each of length $k+(n-k)/s$, which together contain every $(k+1)$-long sequence of distinct elements of $X$ exactly once. This is a \df{narrow} design. Similarly, one could ask for a sequence of length $k + \lambda(n-k)/s$ for some $1\le \lambda\le s$, such that each orbit invariant occurs exactly $\lambda$ times. The corresponding design will contain every element of $X^{(k+1)}$ exactly $\lambda$ times.

As a simple example, the dihedral group $D_n$ of order $2n$ acts $2$-flatly $1$-transitively on $X=\Z/n\Z$ provided $n$ is odd. The orbit invariant on a pair $(x,y)\in X^{(2)}$ is given by $|x-y|$.
Accordingly, the sequencing $(0,1,-1,2,-2,\cdots)$ of the affine action provides a long sequencing for the dihedral action; accordingly, the orbit of this sequencing under $D_n$ is a wide design, containing every ordered pair of distinct elements exactly twice. (In the case $n=5$, this design is exhibited as the rows of the first and fourth latin squares in the example at the end of section \ref{sec:affine} (Table \ref{table:aff5}). Furthermore, fortuitously, taking just the first $k+(n-k)/s = (n+1)/2$ values of this sequence yields a short sequencing since the $(n-1)/2$ absolute values occur exactly once. Accordingly, just looking at the first $(n+1)/2$ entries in the long design provides a narrow design, in which every ordered pair of distinct elements appears exactly once. These are the entries to the left of the dashed lines in Table \ref{table:aff5}.

We remark that the ``designs'' obtained from non-sharply transitive group actions can be thought of as ordered analogues of balanced incomplete block designs. The set $X$ provides the $v=n$ points; the $b=|G|$ (ordered) ``blocks'' are given by the orbit of the sequencing; each (ordered) $t$-subset, with $t=(k+1)$, appears consecutively in exactly $\lambda$ blocks, where $1\le\lambda\le s$, where $s$ is the flatness. These designs seem harder to construct in general than in the usual unordered case, if a suitable group action is not available.

\section{Summary and questions}
We have defined $(n,t)$-permutation designs, and shown their existence for parameters
$(p,3)$ for prime $p$, $(n,n-1)$ for even $n$, and a number of small examples.
However, this paper raises many more questions than it answers.  Here are a few:
\begin{enumerate}
\item For which parameters $(n,t)$ does an $(n,t)$-permutation design exist? An $(n,t)$-\carth rectangle?
\item Study the corresponding covering problem: find the minimal number of permutations of an $n$-set which
   together contain each distinct $t$-tuple as a consecutive subsequence.
\item Construct an explicit family of sequencings of the action of $\Aff{\gf q}$ when $q$ is not prime.
\item Construct an explicit family of sequencings for the action of $\pgl(\gf q)$.
\item Show, for example by strengthening the estimate (\ref{eq:expseq}), that there are only finitely many
  non-sequenceable sharply $k$-transitive actions when $k\ge 2$.
\item Prove Conjecture \ref{conj:existence}, namely that the only non-sequenceable sharply $k$-transitive actions
  with $k\ge 2$ have been identified in this paper.
\item Extend Table \ref{table:dfscounts}. (Remark: the number of sequencings of the odd-order cyclic groups appears in the OEIS \cite[A141599]{OEIS}, but no general formula is given there.)
\item Work out in detail the analogues of distinguishing polynomials for other families of group actions.
\item Construct some $(n,k)$-permutation designs with $k>2$ which do not arise from group actions.
\item As an example of the previous problem, are there situations in which one can construct a permutation
design from a non-sharply $k$-transitive
group action by splitting it? That is, if the pointwise stabilizer of a $k$-tuple has size $s$, when can we
select a fraction $1/s$ of the rows of a long design to get a permutation design?
\end{enumerate}

\bibliographystyle{amsplain}
\bibliography{myRefs}

\appendix
\section{Proof of Theorem \ref{thm:carthage}}
\label{sec:carthageproof}
\proof First suppose $R$ arises from a sharply $t$-transitive action of a group $G$. Let the $2\times(t+1)$
submatrices $A$ and $B$
agree in $t$ columns and in at least one entry of the remaining column; call the corresponding row the "top" row.
By sharp $t$-transitivity, there is a unique $g\in G$ sending the top row of $A$ (or $B$) in the $t$ columns to the bottom row of $A$ in the same columns. Then $g$ must map the entire top row of $A$ (and $B$) to the bottom row; in particular, the remaining element must agree.

Conversely, suppose we have a $t$-\carth rectangle $R$, on a set $X$ of $n$ symbols, for which the
$t$-rectangle condition holds. We can assume the first row of $R$ is just $(1,\ldots,n)$. Identify each row
with the unique permutation in $\Sym{X}$ mapping the first row to it (so the first row corresponds to
$\mathrm{id}_G$.)
We show that if two rows correspond to permutations $\rho$ and $\sigma$ respectively,
then $\sigma\circ\rho^{-1}$ also corresponds to a row of $R$.
 
The first $t$ elements of rows $\rho$ and $\sigma$ are
$(\rho(1),\ldots, \rho(t))$ and $(\sigma(1),\ldots,\sigma(t))$ respectively.
Now, the entries $\rho(1),\ldots,\rho(t)$ appear in the first row in the corresponding columns; by the CR condition,
there is a unique row containing $\sigma(1),\ldots,\sigma(t)$ in those same columns (in the same order).
That row is some permutation $\pi$ of $(1,\ldots,n)$; we claim that $\pi=\sigma\circ\rho^{-1}$. Certainly if
$1\le i \le t$, we have $\pi(\rho(i)) = \sigma(i)$ by construction, so assume $i>t$. Let $A$ be the $2\times(t+1)$
submatrix formed from columns $1,\ldots,t,i$ in rows $\rho$ and $\sigma$. Let $B$ be the $2\times(t+1)$ submatrix
of rows $\mathrm{id}_G$ and $\pi$, choosing the columns so that the first rows of $A$ and $B$ agree. We then have

$$A = \begin{pmatrix}
\rho(1) & \rho(2) & \cdots & \rho(t) & \rho(i)\\
\sigma(1)&\sigma(2)&\cdots & \sigma(t)&\sigma(i)
\end{pmatrix};\quad
B= \begin{pmatrix}
\rho(1) & \rho(2) & \cdots & \rho(t) & \rho(i)\\
\sigma(1)&\sigma(2)&\cdots & \sigma(t)&\pi(\rho(i))
\end{pmatrix}
$$
Then $\sigma(i)=\pi(\rho(i))$ by the $t$-rectangle condition. Hence $\pi=\sigma\circ\rho^{-1}$ corresponds to a
row of $R$. It follows that the set of permutations in $R$ is closed under both inversion and composition, and so
form a group $G$. This group is sharply $t$-transitive by the CR condition.\qed

\section{Explicit sequencings of small sharply transitive actions}
\label{sec:explicitsequencings}
Table \ref{table:sequences} provides an explicit sequencing of every sharply transitive group action listed in the table of primitive permutation group actions in \cite[VI.1.4]{colbourn2010crc}, except for the cases in which such a sequencing does not exist. We have omitted the alternating groups $\Alt n$, which are sequenceable iff $n$ is even, and the symmetric groups $\textrm{Sym}(n)$, which are sequenceable for all $n$. (We have noted that the action of $\pgl(\gf4)$ is not sequenceable; that action doesn't appear in this table since it is the same as $\Alt{\PP^1(\gf4)}$.)

Row $12.4$ of this table exhibits a sequencing of the sharply $5$-transitive action of the
Mathieu group $M_{12}$. The remarkable $(12,6)$-permutation design thus obtained is a list of $95,\!040$
$12$-permutations containing each ordered distinct $6$-tuple exactly once.

\begin{table}
\footnotesize
\begin{tabular}{|c|c|c|p{4in}|}
\hline
ID & $n$ & $k$ & example sequencing\\
\hline
 5.3 & 5 & 2 & \{1, 2, 5, 3, 4\} \\
 6.2 & 6 & 3 & \text{None} \\
 7.4 & 7 & 2 & \{1, 2, 5, 4, 3, 6, 7\} \\
 8.1 & 8 & 2 & \{1, 2, 7, 5, 4, 3, 8, 6\} \\
 8.4 & 8 & 3 & \text{None} \\
 9.3 & 9 & 2 & \{1, 2, 6, 7, 8, 9, 4, 5, 3\} \\
 9.4 & 9 & 2 & \{1, 2, 4, 7, 8, 9, 5, 6, 3\} \\
 10.5 & 10 & 3 & \{1, 2, 3, 9, 5, 4, 8, 6, 7, 10\} \\
 10.6 & 10 & 3 & \{1, 2, 3, 10, 8, 5, 7, 4, 9, 6\} \\
 11.6 & 11 & 4 & \{1, 2, 3, 4, 9, 8, 11, 6, 7, 10, 5\} \\
 12.2 & 12 & 3 & \{1, 2, 3, 9, 12, 4, 5, 8, 6, 7, 11, 10\} \\
 12.4 & 12 & 5 & \{1, 2, 3, 4, 5, 11, 7, 6, 9, 8, 12, 10\} \\
 13.6 & 13 & 2 & \{1, 2, 12, 4, 7, 3, 6, 5, 13, 8, 11, 9, 10\} \\
 14.2 & 14 & 3 & \{1, 2, 3, 11, 10, 8, 9, 6, 4, 5, 14, 12, 13, 7\} \\
 16.3 & 16 & 2 & \{1, 2, 10, 13, 5, 12, 11, 9, 8, 7, 16, 14, 3, 15, 4, 6\} \\
 17.5 & 17 & 2 & \{1, 2, 13, 6, 7, 5, 4, 14, 9, 12, 3, 11, 10, 17, 8, 15, 16\} \\
 17.6 & 17 & 3 & \{1, 2, 3, 15, 7, 14, 11, 13, 17, 10, 4, 16, 6, 9, 5, 12, 8\} \\
 18.2 & 18 & 3 & \{1, 2, 3, 16, 11, 15, 13, 10, 18, 9, 4, 6, 12, 14, 5, 17, 8, 7\} \\
 19.6 & 19 & 2 & \{1, 2, 17, 15, 10, 3, 8, 13, 19, 18, 14, 9, 7, 11, 6, 12, 16, 4, 5\} \\
 20.2 & 20 & 3 & \{1, 2, 3, 19, 7, 10, 8, 16, 14, 12, 5, 13, 18, 6, 17, 4, 11, 20, 9, 15\} \\
 23.4 & 23 & 2 & \{1, 2, 21, 16, 17, 19, 7, 23, 11, 6, 18, 20, 9, 4, 22, 13, 8, 10, 3, 5, 12, 14, 15\} \\
 24.2 & 24 & 3 & \{1, 2, 3, 20, 12, 10, 24, 5, 15, 17, 14, 9, 18, 4, 11, 19, 21, 16, 6, 7, 13, 23, 22, 8\} \\
 25.13 & 25 & 2 & \{1, 2, 14, 9, 16, 23, 12, 20, 21, 13, 11, 4, 15, 18, 22, 3, 10, 25, 17, 5, 19, 8, 24, 6, 7\} \\
 25.14 & 25 & 2 & \{1, 2, 10, 6, 11, 15, 4, 12, 8, 17, 25, 7, 5, 3, 19, 14, 24, 9, 23, 18, 20, 13, 16, 21, 22\} \\
 25.15 & 25 & 2 & \{1, 2, 14, 6, 4, 17, 22, 9, 23, 20, 3, 15, 8, 10, 18, 19, 7, 25, 24, 11, 21, 16, 12, 5, 13\} \\
 26.3 & 26 & 3 & \{1, 2, 3, 4, 19, 9, 6, 22, 12, 16, 10, 24, 25, 21, 13, 20, 7, 18, 5, 8, 26, 17, 14, 11, 15, 23\} \\
 26.4 & 26 & 3 & \{1, 2, 3, 18, 4, 14, 6, 25, 24, 22, 13, 7, 19, 23, 9, 17, 26, 12, 16, 8, 15, 21, 11, 10, 20, 5\} \\
 27.6 & 27 & 2 & \{1, 2, 24, 17, 23, 11, 13, 5, 20, 26, 22, 16, 18, 8, 6, 12, 21, 3, 14, 25, 7, 9, 27, 10, 4, 15, 19\} \\
 28.7 & 28 & 3 & \{1, 2, 3, 14, 6, 8, 25, 19, 21, 20, 18, 17, 7, 23, 24, 27, 13, 22, 9, 5, 15, 28, 10, 11, 4, 12, 26, 16\} \\
 29.1 & 29 & 2 & \{1, 2, 23, 13, 7, 8, 25, 4, 21, 11, 28, 20, 9, 16, 6, 22, 5, 24, 29, 17, 12, 10, 27, 15, 26, 3, 14, 18, 19\} \\
 30.2 & 30 & 3 & \{1, 2, 3, 10, 9, 26, 13, 21, 14, 5, 7, 19, 20, 8, 30, 11, 23, 24, 4, 22, 15, 18, 16, 6, 17, 27, 12, 25, 29, 28\} \\
 31.8 & 31 & 2 & \{1, 2, 7, 19, 15, 21, 6, 16, 3, 13, 11, 28, 14, 8, 31, 25, 12, 24, 4, 26, 5, 20, 23, 9, 27, 30, 22, 29, 10, 17, 18\}
   \\
 32.1 & 32 & 2 &
   \{1, 2, 22, 20, 31, 11, 28, 3, 5, 21, 27, 10, 8, 16, 24, 13, 14, 9, 15, 18, 17, 6, 32, 7, 23, 12, 25, 30, 26, 19, 4, 29\} \\
 32.4 & 32 & 3 &
   \{1, 2, 3, 28, 25, 15, 27, 7, 21, 10, 24, 13, 14, 26, 11, 20, 23, 32, 4, 29, 5, 12, 17, 19, 16, 18, 8, 9, 31, 30, 6, 22\} \\
 33.1 & 33 & 3 &
   \{1, 2, 3, 12, 7, 29, 19, 10, 5, 31, 8, 22, 21, 9, 26, 13, 14, 17, 27, 24, 23, 4, 11, 20, 25, 18, 15, 30, 16, 6, 33, 32, 28\}\\
  37.1 & 37 & 2 &
   \{1, 2, 9, 13, 3, 33, 25, 4, 12, 7, 35, 28, 19, 16, 21, 15, 26, 23, 20, 27, 29, 32, 24, 34, 17, 8, 37, 36, 18, 6, 10, 11, 22, 14, 5, 30, 31\}\\
  38.2 & 38 & 3 &
   \{1, 2, 3, 11, 27, 29, 14, 30, 6, 15, 22, 13, 38, 9, 24, 32, 21, 31, 37, 8, 23, 25, 26, 4, 17, 35, 36, 7, 16, 10, 12, 5, 20, 28, 19, 18, 34, 33\} \\
 \hline
\end{tabular}
\caption{\strut Sequencings of the sharply transitive group actions listed in \cite[VI.1.4]{colbourn2010crc}.}
\label{table:sequences}
\end{table}

\begin{table}
\begin{center}
\footnotesize
\begin{tabular}{|c|c|c|p{4in}|}
\hline
ID & $n$ & $k$ & example sequencing\\
 \hline
 41.8 & 41 & 2 &
   \{1, 2, 5, 6, 10, 16, 20, 19, 25, 36, 31, 34, 29, 7, 9, 28, 23, 32, 8, 39, 41, 4, 37, 17, 15, 38, 30, 33, 12, 22, 26, 24, 27, 3, 11, 21, 40, 18, 35, 13, 14\} \\
 42.2 & 42 & 3 &
   \{1, 2, 3, 17, 38, 20, 28, 18, 33, 12, 8, 19, 13, 30, 34, 40, 23, 26, 22, 31, 42, 15, 39, 14, 27, 41, 16, 11, 6,  9, 5, 4, 21, 10, 7, 32, 24, 29, 25, 36, 35, 37\} \\
 43.8 & 43 & 2 &
   \{1, 2, 31, 36, 37, 29, 19, 39, 30, 40, 11, 10, 15, 26, 42, 14, 9, 23, 12, 21, 33, 22, 13, 4, 5, 18, 20, 41, 6, 3, 17, 43, 16, 24, 35, 27, 38, 25, 32, 34, 28, 7, 8\} \\
 44.2 & 44 & 3 &
   \{1, 2, 3, 30, 35, 19, 40, 21, 38, 13, 23, 5, 36, 6, 17, 43, 41, 9, 22, 27, 42, 11, 28, 18, 14, 44, 33, 4, 29, 12, 10, 26, 7, 39, 16, 31, 20, 24, 8, 37, 34, 32, 25, 15\} \\
 47.4 & 47 & 2 &
   \{1, 2, 4, 20, 40, 43, 6, 13, 38, 24, 27, 9, 29, 5, 10, 15, 17, 26, 46, 11, 31, 28, 14, 32, 33, 44, 34, 37, 36, 22, 18, 45, 19, 35, 39, 41, 30, 12, 47, 3, 16, 23, 21, 25, 42, 7, 8\} \\
 48.2 & 48 & 3 &
   \{1, 2, 3, 22, 47, 9, 41, 45, 42, 10, 48, 15, 16, 39, 6, 44, 35, 29, 4, 24, 34, 27, 26, 7, 18, 43, 13, 33, 8, 17, 36, 11, 21, 19, 37, 38, 40, 46, 23, 32, 20, 30, 25, 31, 5, 14, 28, 12\} \\
 49.22 & 49 & 2 &
   \{1, 2, 41, 8, 14, 19, 49, 23, 24, 17, 21, 22, 9, 31, 6, 38, 25, 15, 10, 4, 18, 36, 13, 5, 44, 45, 48, 32, 7, 43, 26, 20, 40, 34, 27, 29, 37, 28, 46, 33, 11, 3, 47, 39, 35, 42, 12, 30, 16\} \\
 49.23 & 49 & 2 &
   \{1, 2, 18, 9, 27, 13, 42, 16, 26, 5, 3, 32, 33, 38, 23, 43, 10, 11, 12, 36, 22, 4, 34, 35, 30, 14, 49, 31, 8, 39, 28, 17, 45, 15, 24, 44, 46, 19, 37, 47, 41, 25, 21, 40, 20, 7, 29, 6, 48\} \\
 49.24 & 49 & 2 &
   \{1, 2, 39, 22, 48, 21, 40, 42, 33, 37, 43, 16, 19, 14, 31, 18, 10, 15, 29, 36, 27, 44, 41, 24, 34, 8, 30, 49, 26, 7, 35, 28, 6, 9, 3, 11, 46, 32, 38, 25, 23, 13, 20, 47, 17, 45, 12, 4, 5\} \\
 50.2 & 50 & 3 &
   \{1, 2, 3, 32, 5, 26, 34, 38, 43, 48, 21, 18, 33, 47, 42, 35, 7, 30, 12, 25, 46, 14, 17, 19, 13, 50, 49, 39, 37, 4, 15, 29, 16, 45, 20, 11, 40, 41, 28, 6, 31, 9, 10, 24, 27, 22, 44, 8, 23, 36\} \\
 50.3 & 50 & 3 &
   \{1, 2, 3, 35, 28, 36, 12, 6, 48, 7, 31, 15, 13, 4, 19, 42, 8, 17, 11, 26, 22, 21, 47, 39, 27, 33, 5, 25, 38, 29, 10, 30, 20, 32, 43, 41, 49, 45, 34, 9, 44, 16, 14, 40, 18, 23, 37, 50, 46, 24\} \\
   \hline
\end{tabular}

\normalsize\strut Table \ref{table:sequences}, cont.
\end{center}
\end{table}

\end{document}